\documentclass{article} [oneside]
\usepackage{array}
\usepackage{amssymb}
\usepackage{amsthm}
\usepackage{xcolor}
\usepackage[leqno]{amsmath}
\usepackage{amsfonts}
\usepackage{graphicx}
\usepackage[toc,page]{appendix}
\usepackage[margin=0pt,font+=small,labelformat=parens,labelsep=space, skip=6pt,list=false,hypcap=false] {caption}
\usepackage{subfig}
\usepackage{relsize}
\usepackage[left]{lineno}
\usepackage{amsthm}

\newtheorem{theorem}{Theorem}[section]
\newtheorem{Lemma}[theorem]{Lemma}
\usepackage[top=.9in,bottom=.9in]{geometry}
 \numberwithin{equation}{section}
 
\begin{document}

\pagestyle{plain}
%\pagenumbering{arabic}
\numberwithin{equation}{section}
%%%
\iftrue
\title{A Curious Trigonometric Infinite Product in Context}
\author{ Michael Milgram\footnote{mike@geometrics-unlimited.com}\\{Consulting Physicist, Geometrics Unlimited, Ltd.}}
%{{Box 1484, Deep River, Ont. Canada. K0J 1P0}}
\maketitle
\begin{flushleft} \vskip 0.3 in % must FOLLOW maketitle
%\vskip .2in
%\centerline{ Michael Milgram\footnote{mike@geometrics-unlimited.com}}
%\centerline{Consulting Physicist, Geometrics Unlimited, Ltd.}
\centerline{Box 1484, Deep River, Ont. Canada. K0J 1P0}
\vskip .2in
\centerline{Author's manuscript, Orcid:0000-0002-7987-0820}
revised \today
\vskip .2in
\fi
\centerline{}
\vskip .1in
MSC classes: 	40A20, 26-02, 26A09, 26E99, 33-02, 33B10, 
\vskip 0.1in
Keywords:
Infinite product, Vi\`ete's Product, Accumulation point, Weierstrass' Factor theorem, trigonometric products
\vskip 0.1in 

\centerline{\bf Abstract}\vskip .3in

By treating the multiple argument identity of the logarithm of the Gamma function as a functional equation, we obtain a curious infinite product representation of the $sinc$ function in terms of the cotangent function. This result is believed to be new. It is then shown how to convert the infinite product to a finite product, which turns out to be a simple telescoping of the double angle $sin$ function. In general, this result unifies known infinite product identities involving various trigonometric functions when the product term index appears as an exponent. In one unusual case, what appears to be a straightforward limit, suggests a counterexample to Weierstrass' factor theorem. A resolution is offered. An Appendix presents the general solution to a simple functional equation. This work is motivated by its educational interest.

\section{Introduction and Motivation}

% from file Plimits2.mw in directory Mellints
Throughout, the symbols, $m,n,j,k,N$ are positive integers and $\gamma$ is the Euler-Mascheroni constant; other symbols are real, although the extension to general complex values is always available. In the course of a study of the properties of $\eta(s)$, the alternating companion to Riemann's function $\zeta(s)$ defined \cite[Eq. 2.3(1)]{Sriv&Choi} by

\begin{equation}
\eta(s)\equiv\sum_{k=0}^{\infty}(-1)^{k}/(k+1)^s=(1-2^{1-s})\,\zeta(s),\hspace{15pt} \Re(s)>0\,,
\label{EtaDef}
\end{equation} 

I obtained the identity

\begin{equation}
\overset{\infty}{\underset{j =1}{\sum}}\; \frac{\eta \! \left(1+j \right) \left(-a \right)^{j}}{1+j}
 = 
\frac{1}{a}\ln \! \left(\frac{\sqrt{\pi}\,\Gamma \left(\frac{a}{2}+1\right)}{\Gamma \left(\frac{a}{2}+\frac{1}{2}\right) 2^{a}}\right) ,\hspace{10Pt} |a|<1
\label{Cm1a}
\end{equation}
corresponding to similar identities developed in \cite[Section 3.4]{Sriv&Choi}, in particular the identity \cite[Eq. 3.4(10)]{Sriv&Choi} 

\begin{equation}
\overset{\infty}{\underset{j =1}{\sum}}\; \frac{\zeta \! \left(1+j \right) \left(-a \right)^{j}}{1+j}
 = -\frac{\ln \! \left(\Gamma \! \left(a +1\right)\right)}{a}-\gamma,\hspace{10Pt} |a|<1\,.
\label{Sc1a}
\end{equation}

To obtain \eqref{Cm1a}, apply the second equality of \eqref{EtaDef} to yield
\begin{equation}
\overset{\infty}{\underset{j =1}{\sum}}\; \frac{\zeta \! \left(1+j \right) \left(-a \right)^{j}}{1+j}
 = 
\frac{1}{a}\ln \! \left(\frac{\sqrt{\pi}\,\Gamma \left(\frac{a}{2}+1\right)}{\Gamma \left(\frac{a}{2}+\frac{1}{2}\right) 2^{a}}\right)+\overset{\infty}{\underset{j =1}{\sum}}\; \frac{\zeta \! \left(1+j \right) \left(-\frac{a}{2}\right)^{j}}{1+j}
\label{Cm1B}
\end{equation}

from which \eqref{Cm1a} can be easily verified by recourse to \eqref{Sc1a}. In addition to employing \eqref{Cm1B} to prove \eqref{Cm1a}, cursory observation of \eqref{Cm1B} shows that it satisfies a simple and easily solvable (see Appendix) functional equation of the form 

\begin{equation}
f \! \left(a \right) = 
g \! \left(a \right)+x\,f \! \left(\frac{a}{p}\right)\,,
\label{R0a}
\end{equation}
of which the application to $\eta(s)$ is discussed elsewhere \cite{Mprep}. However, noticing that 
\eqref{Sc1a} and \eqref{Cm1B} also simply reproduce the half argument identity of $\log(\Gamma(1+a))$ in the form
\begin{equation}
\ln \! \left(\Gamma \! \left(1+a \right)\right) = 
\ln \! \left(2\right) a -\frac{\ln \! \left(\pi \right)}{2}+\ln \! \left(\Gamma \! \left(\frac{a}{2}+\frac{1}{2}\right)\right)+\ln \! \left(\Gamma \! \left(1+\frac{a}{2}\right)\right)\,,
\label{DargId}
\end{equation}
the question arises: ``Will anything interesting be found if \eqref{DargId} is studied as an instance of the functional equation \eqref{R0a}?" That is the motivation for this work, which is intended to be an educational exercise in the merits of curiosity.\newline

In the following Section \ref{sec:FeqSol}, that question is pursued, first in Subsection \ref{sec:Curious} where it is shown to yield interesting and curious infinite trigonometric products and corresponding infinite summation identities, some of which are thought to be new. In that same Subsection these identities are related to known results by the operation of differentiation and the exceptional cases are analysed. Curiosity then suggests that the new, curious identity be generalized, and this is done in two different ways in Sections \ref{sec:genRecur} and \ref{sec:TripOn}; in the latter Section, connections are also made with known identities, thereby unifying known  results listed elsewhere. A requirement that always arises when studying what may be a new identity is to place it in context with known material from the literature; this is pursued in Section \ref{sec:Perspective} where a closely-related, purported, and long-ago analysed, identity is reviewed. Finally, curiosity leads one far afield by noticing that the limiting case of two of the contextually related identities previously introduced ostensibly appear to violate Weierstrass' factor theorem. In Section \ref{sec:Digression} this issue is considered and a resolution is offered.

\section{$\log (\Gamma(a+1))$ treated as a functional equation}\label{sec:FeqSol}

\subsection{A curious identity} \label{sec:Curious}
By applying the first corollary of Lemma \eqref{R0} to \eqref{DargId}-- see the Appendix-- we easily find the interesting identity

\begin{equation}
\Gamma \! \left(1+2\,a \right)\overset{\infty}{\underset{j =0}{\prod}}\; \frac{1}{2^{\,2\,a}}\left(\frac{\sqrt{\pi}\,\Gamma \left({a}/{2^{j}}+1\right)}{\Gamma \left({a}/{2^{j}}+\frac{1}{2}\right)}\right)^{2^{ j}}={{\mathrm e}^{-2\,\gamma \,a}}\,,
\label{R0A}
\end{equation}
equivalent to the (convergent) sum
\begin{equation}
\overset{\infty}{\underset{j =1}{\sum}}\; \left(2^{j}\,\ln \! \left(\frac{\Gamma \! \left({a}/{2^{j}}+1\right)}{\Gamma \! \left({a}/{2^{j}}+{1}/{2}\right)}\right)+2^{j -1}\,\ln \! \left(\pi \right)-2\,\ln \! \left(2\right) a \right)+2\,\ln \! \left(\Gamma \! \left(a +1\right)\right)
 = -2\,\gamma \,a\,.
\label{Rs2}
\end{equation}
% From file ExId.mw

Now add the identity \eqref{Rs2} and its companion where $a\Rightarrow -a$, to obtain the curious infinite product
\begin{equation}
\overset{\infty}{\underset{j =1}{\prod}}\! \left(\tan \! \left({a}/{2^{j}}\right)\frac{2^{j}\,}{a}\right)^{2^{j -1}}
 = \frac{a}{\sin \! \left(a \right)}\,,
\label{Vsum2}
\end{equation}
and, letting $a=ib$, its hyperbolic equivalent
\begin{equation}
\overset{\infty}{\underset{j =1}{\prod}}\; \left(\frac{2^{j}\,\tanh \! \left({b}/{2^{j}}\right)}{b}\right)^{2^{j -1}}
 = \frac{b}{\sinh \! \left(b \right)}\,.
\label{Vsum2a}
\end{equation}
Inverting also yields the following representation of the $\mathrm{sinc}(a)$ function, defined by $\mathrm{sinc} (a) \equiv \sin(a)/a$ :
\begin{equation}
\mathrm{sinc} \! \left(a \right) = 
\overset{\infty}{\underset{j =1}{\prod}}\; \left(\cot \! \left({a}/{2^{j}}\right)\frac{ a}{2^{j}}\right)^{2^{j -1}}\,.
\label{sinc}
\end{equation}

It is important to consider the exceptional cases where $a$ is in integral multiple of $\pi$. There are two possibilities: $a$ is even or odd. In the latter case, let $a=(2\,n-1)\,\pi,~n\geq 1$, in which case the only divergent term in the product is labelled by the index $j=1$. So, by setting $n\rightarrow n+\epsilon$, consider the limit

\begin{equation}
\underset{\mathit{\epsilon} \rightarrow 0}{\mathrm{lim}}\! \frac{2\,\tan \! \left(\frac{\left(2\,n +2\,\mathit{\epsilon} -1\right) \pi}{2}\right)}{\left(2\,n +2\,\mathit{\epsilon} -1\right) \pi}
 = 
-\frac{2}{\pi^{2} \left(2\,n -1\right) \mathit{\epsilon}}+O\! \left(1\right)
\label{Cp1}
\end{equation}
and compare to the corresponding divergent term on the right-hand side of \eqref{Vsum2}, to obtain
\begin{equation}
\overset{\infty}{\underset{j =2}{\prod}}\; \left(\frac{\tan \! \left(\frac{\left(2\,n -1\right) \pi}{2^{j}}\right) 2^{j}}{\left(2\,n -1\right) \pi}\right)^{2^{j -1}}
 = \pi^{2} \left(1/2-n \right)^{2}\,.
\label{CpOdd}
\end{equation}
In the case that $a$ is even, recall that any even number can always be written as the product of an integral power of $2$ and an odd number, so let $a=2^m\,(2n-1)$, where $m\geq 0$ and $n\geq 1$. By setting $m\rightarrow m+\epsilon$ and considering the limit as $\epsilon\rightarrow 0$, it is easy to find that the first $m$ terms of the product each vanish to order $(\epsilon\,\ln(2))^{2^{j-1}}$, and therefore
\begin{equation}
\underset{\epsilon \rightarrow 0}{\mathrm{lim}}\, \;\overset{m}{\underset{j =1}{\prod}}\; \left(\frac{\tan \! \left(\frac{2^{m +\epsilon}\,(2n-1)\,\pi}{2^{j}}\right) 2^{j}}{2^{m +\epsilon}\,(2n-1)\,\pi}\right)^{2^{j -1}}
\approx \frac{\epsilon^{2^{m}}\,\ln \! \left(2\right)^{2^{m}-1}}{\epsilon}\,.
\label{Plim1}
\end{equation}
The next term of the product corresponding to $j=m+1$ is the only divergent factor, being of order 
\begin{equation}
\underset{\epsilon \rightarrow 0}{\mathrm{lim}}\; \left(\frac{\tan \! \left(\frac{2^{m +\epsilon}\,(2n-1)\,\pi}{2^{m +1}}\right) 2^{m +1}}{2^{m +\epsilon}\,(2n-1)\,\pi}\right)^{2^{m}}
 \approx 
\left(\frac{4}{\ln \! \left(2\right) (2n-1)^{2}\,\pi^{2}\,\epsilon}\right)^{2^{m}}\,.
\label{Pe01a}
\end{equation}
Therefore the product of \eqref{Plim1} and \eqref{Pe01a} diverges to order $\epsilon^{-1}$, as does the right-hand side of \eqref{Vsum2} under the conditions specified. So, if $a$ is an integral multiple of $\pi$, we have the exceptional case

\begin{equation}
\overset{\infty}{\underset{j =m +2}{\prod}}\; \left(\frac{\tan \! \left(2^{m -j}\,(2n-1)\,\pi \right)}{2^{m -j}\,(2n-1)\,\pi}\right)^{2^{j -1}}
 = \left(\frac{(2n-1)\,\pi}{2}\right)^{2^{m +1}}
\label{PEO2}
\end{equation}
reducing to \eqref{CpOdd} if $m=0$.\newline

 {\bf Remarks:}
 \begin{itemize}
 \item Differentiating \eqref{Vsum2} with respect to $a$ gives
\begin{equation}
\overset{\infty}{\underset{j =1}{\sum}}\! \left(2^{j -1}-\frac{a}{\sin \! \left(a/2^{j-1} \right)}\right)
 = a\,\cot \! \left(a \right)-1\,,
\label{Vsum3}
\end{equation}
an extension of the listed \cite[Eq. 25.1.2]{Hansen} identity
\begin{equation} 
\overset{n}{\underset{j =0}{\sum}}\; \frac{1}{\sin \! \left({x}/{2^{j}}\right)}
 = \cot \! \left(\frac{x}{2^{n +1}}\right)-\cot \! \left(x \right)\,.
\label{H25p1p2}
\end{equation}

 \item It is possible to reproduce \eqref{Vsum3} by subtracting the sum \eqref{H25p1p2} from the elementary identity
 \begin{equation}
 \overset{n}{\underset{j =0}{\sum}}\; 2^{j} = 2^{n +1}-1
 \label{Also}
 \end{equation}
 and evaluating the limit as $n\rightarrow\infty$. 
 
 \item Further, differentiation of \eqref{Rs2} yields
\begin{equation}
\overset{\infty}{\underset{j =1}{\sum}}\! \left(2\,\ln \! \left(2\right)-\psi \! \left(1+\frac{a}{2^{j}}\right)+\psi \! \left(\frac{a}{2^{j}}+\frac{1}{2}\right)\right)
 = 2\,\psi \! \left(a +1\right)+2\,\gamma\,,
\label{R1bD}
\end{equation} 
a result that is discussed elsewhere \cite{Mprep}, which could also be obtained by evaluating the $n\rightarrow \infty$ limit of \cite[Eq. (55.8.4)]{Hansen}. 

\item A related result extracted from the classical literature \cite[Eq. (1042)]{Jolley} is
\begin{equation}
\underset{n \rightarrow \infty}{\mathrm{lim}}\! \left(\overset{n}{\underset{j =1}{\prod}}\; \tan \! \left(\frac{\pi \,j}{2\,n}\right)\right)^{\frac{1}{n}}
 = 1\,.
\label{Jo1}
\end{equation}

\end{itemize}

\subsection{Generalizations} \label{sec:gen}

\subsubsection{By recursion} \label{sec:genRecur}

% from file Trigproducts.mw
It is possible to generalize \eqref{Vsum2} as follows
\begin{equation}
\overset{\infty}{\underset{j =n +1}{\prod}}\; \left(\frac{2^{j-n}\,}{a}\tan \! \left({2^{n-j}\,a}\right)\right)^{2^{j -1}}
 = \frac{a^{2^{n}}}{\sin^{2^{n}}\left(a \right)},\hspace{10pt} n\geq 0,
\label{Gp1b}
\end{equation}
the proof being by induction.\newline

{\bf Proof:} Let $n\Rightarrow n+1$ and simultaneously let $a\Rightarrow a/2$ to obtain
\begin{equation}
\overset{\infty}{\underset{j =n +2}{\prod}}\; \left(\frac{2^{j-n}}{a}\,\tan \! \left(2^{n -j}\,a \right)\right)^{2^{j -1}}
 = 
\frac{\left(\frac{a}{2}\right)^{2^{n+1}}}{\sin^{2^{n+1}}\left(\frac{a}{2}\right)}.
\label{P1}
\end{equation}
Extend the product on the left-hand side and transfer the new factor to the right-hand side, yielding
\begin{equation}
\overset{\infty}{\underset{j =n +1}{\prod}}\; \left(\frac{2^{j-n }}{a}\,\tan \! \left(2^{n -j}\,a \right)\right)^{2^{j -1}}
 = 
\frac{\left(\frac{a}{2}\right)^{2^{n+1}}}{\sin^{2^{n+1}}\left(\frac{a}{2}\right)} \left(\frac{2\,\tan \left(\frac{a}{2}\right)}{a}\right)^{2^{n}}\,.
\label{P2}
\end{equation}
After some elementary simplification we find

\begin{equation}
\overset{\infty}{\underset{j =n +1}{\prod}}\; \left(\frac{2^{j-n}}{a}\,\tan \! \left(2^{n -j}\,a \right)\right)^{2^{j -1}}
 = 
a^{2^{n}} \left(\frac{1}{2\,\cos \! \left(\frac{a}{2}\right) \sin \! \left(\frac{a}{2}\right)}\right)^{2^{n}}
\label{P3}
\end{equation}
reducing to \eqref{Gp1b} from the double angle identity 
\begin{flalign}
\sin(a)=2\cos(a/2)\sin(a/2). \hspace{50pt} \mbox{ \bf{QED} }
\label{sin2a}
\end{flalign}

Again generalizing, let $n\Rightarrow n_{1}, a\Rightarrow a/2^{n_{1}}$ in \eqref{Gp1b} and similarly $n\Rightarrow n_{2}, a\Rightarrow a/2^{n_{2}}$ with $n_{2}>n_{1}$. Divide the resulting identities and split the first product from $n_{1}$ to $n_{2}$ and $n_{2}+1$ to $\infty$. After some simplification, we have the finite product

\begin{equation}
\overset{n_{2}}{\underset{j =n_{1}+1}{\prod}}\; \left(\frac{2^{j}}{a}\,\tan \! \left(\frac{a}{2^{j}}\right)\right)^{2^{j -1}}
 = 
a^{N_{1}-N_{2}}\,2^{-N_{1}\,n_{1}+N_{2}\,n_{2}}\frac{ \sin^{N_{2}}\!\left(\frac{a}{N_{2}}\right)}{\sin^{N_{1}}\!\left(\frac{a}{N_{1}}\right)}\,,
\label{P5}
\end{equation}
where $N_{1}\equiv 2^{n_{1}}$ and  $N_{2}\equiv 2^{n_{2}}$. Consider two examples. Let $n_{1}=0,~n_{2}=4$ and $a\Rightarrow 2\,a$ giving
\begin{equation}
\tan \! \left(a \right) \tan^{2}\left(\frac{a}{2}\right)\tan^{4}\left(\frac{a}{4}\right)\tan^{8}\left(\frac{a}{8}\right)
 = 
 2^{15}\frac{\sin^{16}\left(\frac{a}{8}\right)}{\sin \! \left(2\,a \right)}\,.
\label{X1}
\end{equation}
The case $n_{1}=n,~n_{2}=n+1$ yields the identity
\begin{equation}
\tan \! \left(\frac{a}{2^{n +1}}\right) = 
\frac{2 \sin^{2}\left({a}/{2^{n +1}}\right)}{\sin \! \left({a}/{2^{n}}\right)}\,,
\label{X1b}
\end{equation}
made obvious by setting $a\Rightarrow 2^n\,a$ and referring to \eqref{sin2a}. It is easily seen that \eqref{X1} is a simple recursive rendition of the double angle formula \eqref{sin2a} by rewriting it in the form

\begin{equation}
\sin \! \left(2\,a \right)\frac{ \sin \! \left(a \right)}{\cos \! \left(a \right)}\tan^{2}\left(\frac{a}{2}\right) \tan^{4}\left(\frac{a}{4}\right)\tan^{8}\left(\frac{a}{8}\right)
 =  2^{15}\sin^{16}\left(\frac{a}{8}\right)
\label{X1a}
\end{equation}
and applying \eqref{sin2a} recursively.

\subsubsection{Triplication and Onwards} \label{sec:TripOn}
% from file TrigProducts.mw

Since \eqref{R0A} is essentially based on the duplication formula \eqref{DargId}, consider the triplication analogue
\begin{equation}
\ln \! \left(\Gamma \! \left(a +1\right)\right) = 
\ln \! \left(\frac{3^{a +\frac{1}{2}}\,\Gamma \! \left(\frac{a}{3}+\frac{1}{3}\right) \Gamma \! \left(\frac{a}{3}+\frac{2}{3}\right)}{2\,\pi}\right)+\ln \! \left(\Gamma \! \left(\frac{a}{3}+1\right)\right)
\label{Gn3}
\end{equation}

a special case of the well-known general identity \cite[Eq. 5.5.6]{NIST}
\begin{equation}
\Gamma \! \left(n\,a \right) = 
\left(2\,\pi \right)^{\frac{1}{2}-\frac{n}{2}}\,n^{na -\frac{1}{2}} \,\overset{n-1}{\underset{k =0}{\prod}}\; \Gamma \! \left(a +\frac{k}{n}\right)\,.
\label{Gn}
\end{equation}
Apply Corollary 1 from the Appendix as before to find

\begin{equation}
\ln \! \left(\Gamma \! \left(a +1\right)\right) = 
\overset{\infty}{\underset{j =0}{\sum}}\; \left(\frac{\ln \! \left(3\right)}{2}+\ln \! \left(3^{^{\frac{3\,a}{3^{1+j}}}}\right)+\ln \! \left(\Gamma \! \left(\frac{a}{3^{1+j}}+\frac{1}{3}\right)\right)+\ln \! \left(\Gamma \! \left(\frac{a}{3^{1+j}}+\frac{2}{3}\right)\right)-\ln \! \left(2\,\pi \right)\right)\,,
\label{Gn3a}
\end{equation}
the analogue of \eqref{Rs2}. Adding \eqref{Gn3} and its partner with $a\Rightarrow -a$ then yields
\begin{equation}
\overset{\infty}{\underset{j =0}{\prod}}\; \frac{1+2\,\cos \! \left({2\,a}/{3^{1+j}}\right)}{3}=\frac{\sin \! \left(a \right)}{a} \,,
\label{Gn3cA}
\end{equation}
equivalent to
\begin{equation}
\overset{\infty}{\underset{j =1}{\prod}}\! \left(1-\frac{4 \sin^{2}\left({a}/{3^{j}}\right)}{3}\right)=\frac{\sin \! \left(a \right)}{a}\,,
\label{Gn3cI}
\end{equation}
a known result \cite[Eq. 91.8.14]{Hansen}\,. Similarly, by differentiating \eqref{Gn3a} we find
\begin{equation}
\overset{\infty}{\underset{j =0}{\sum}}\; \left(\frac{3\,\ln \! \left(3\right)}{3^{1+j}}+\frac{\psi \! \left({a}/{3^{1+j}}+\frac{1}{3}\right)}{3^{1+j}}+\frac{\psi \! \left({a}/{3^{1+j}}+\frac{2}{3}\right)}{3^{1+j}}\right)
 = \psi \! \left(a +1\right)\,,
\label{Gn3aD}
\end{equation}
 the analogue of \eqref{R1bD}. {\bf Remark:} Hansen \cite{Hansen} sources \eqref{Gn3cI} by the initials ``CJ', which label does not appear in the Bibliography.\newline
 
Continuing, and similarly, the quadruplication and quintuplication cases produce 

\begin{equation}
\overset{\infty}{\underset{j =0}{\prod}}\; \left(\frac{\cos \! \left({a}/{4^{1+j}}\right)}{2}+\frac{\cos \! \left({3\,a}/{4^{1+j}}\right)}{2}\right)=\frac{\sin \! \left(a \right)}{a}
\label{Gn4A}
\end{equation}
\iffalse 
\begin{equation}
\overset{\infty}{\underset{j =0}{\sum}}\! \left(\frac{8\,\ln \! \left(2\right)}{4^{1+j}}+\frac{\psi \! \left({a}/{4^{1+j}}+\frac{1}{4}\right)}{4^{1+j}}+\frac{\psi \! \left({a}/{4^{1+j}}+\frac{1}{2}\right)}{4^{1+j}}+\frac{\psi \! \left({a}/{4^{1+j}}+\frac{3}{4}\right)}{4^{1+j}}\right)
 = \psi \! \left(a +1\right)
\label{Gn4b}
\end{equation}
\fi 
and
\begin{equation}
\overset{\infty}{\underset{j =0}{\prod}}\; \frac{1+2\,\cos \! \left({2\,a}/{5^{1+j}}\right)+2\,\cos \! \left({4\,a}/{5^{1+j}}\right)}{5}=\frac{\sin \! \left(a \right)}{a} \,,
\label{Gn5B}
\end{equation}
the general even and odd cases being respectively
\begin{equation}
\overset{\infty}{\underset{j =0}{\prod}}\; \frac{1}{N}\;\overset{N}{\underset{n =1}{\sum}}\;{\cos \! \left(\frac{\left(2\,n -1\right) a}{\left(2\,N \right)^{1+j}}\right)}{}=\frac{\sin \! \left(a \right)}{a} \,,
\label{G2A}
\end{equation}
\begin{equation}
\overset{\infty}{\underset{j =0}{\prod}}\; \left(\frac{1}{2\,N +1}+\frac{2 }{2\,N +1}\;\overset{N}{\underset{n =1}{\sum}}\; \cos \! \left(\frac{2\,n\,a}{\left(2\,N +1\right)^{1+j}}\right)\right)=\frac{\sin \! \left(a \right)}{a} \,.
\label{G2x}
\end{equation}
In \eqref{G2A}, let $N=1$ to find Euler's well-known infinite product formula \cite[Eq. 91.8.2]{Hansen}

\begin{equation}
\overset{\infty}{\underset{j =0}{\prod}}\; \cos \! \left({a}/{2^{1+j}}\right)=\frac{\sin \! \left(a \right)}{a} \,.
\label{GN1}
\end{equation}
Both \eqref{G2A} and \eqref{G2x} have been obtained elsewhere through the use of both Chebyshev polynomials and Fourier transforms (see Nishimura \cite{Nishi2016} where the words $odd$ and $even$ are erroneously scrambled). \newline

{\bf Remarks} 
\begin{itemize}
\item {Although the identity \eqref{GN1} as cited requires $|a|<1$, the derivation given here removes that limitation.}

\item
Additionally, \eqref{GN1} might be compared to the equally well-known \cite[Eq. 91.8.3, unsourced]{Hansen} and occasionally re-cited (e.g. \cite{WikiTrig}) identity
\begin{equation}
\overset{k -1}{\underset{j =0}{\prod}}\; \cos \! \left(2^{j}\,a \right)
 =\, 
\frac{\sin \! \left(2^{k}\,a \right)}{2^{k}\,\sin \! \left(a \right)}\,.
\label{BR114}
\end{equation}
\item
We also find from the classical literature \cite[Eq. (1050)]{Jolley} the related identity
\begin{equation}
\overset{2\,k -1}{\underset{j =1}{\prod}}\; \cos \! \left(\frac{\pi \,j}{k}\right)
 = \frac{\left(-1\right)^{k}-1}{2^{2\,k -1}}\,.
\label{Jo2}
\end{equation}
\end{itemize}

Superficially, the equality of the left-hand sides of \eqref{Gn4A} and \eqref{Gn5B} appears to be surprising; likewise the corresponding sides of \eqref{G2A} and \eqref{G2x}. However, from \cite[Eq. 17.1.1]{Hansen} we have

\begin{equation}
\overset{N}{\underset{n =1}{\sum}}\! \cos \! \left(\frac{\left(2\,n -1\right) a}{\left(2\,N \right)^{1+j}}\right)
 = 
\frac{\sin \! \left({a}/{\left(2\,N \right)^{j}}\right)}{2\,\sin \! \left({a}/{\left(2\,N \right)^{1+j}}\right)}
\label{SumId1}
\end{equation}
and with some simplification help from Maple \cite{Maple}
\begin{equation}
\overset{N}{\underset{n =1}{\sum}}\; \cos \! \left(\frac{2\,n\,a}{\left(2\,N +1\right)^{1+j}}\right)
 = 
\frac{\sin \! \left(a /\left(2\,N +1\right)^{j}\right)}{2\,\sin \! \left(a/ \left(2\,N +1\right)^{1+j}\right)}-\frac{1}{2}\,,
\label{SumId2}
\end{equation} 
so that, after the application of \eqref{SumId2} to \eqref{G2x}, we have
\iffalse
\begin{equation}
\frac{\textcolor{blue}{\sin \! \left(a \right)}\; \textcolor{red}{\sin \! \left(\frac{a}{2\,N\, +1}\right)} \overset{\infty}{\underset{j =2}{\prod}}\;\textcolor{black}{ \frac{\,\sin \left(a/ \left(2\,N\, +1\right)^{j}\right)}{\sin \left(a/ \left(2\,N +1\right)^{1+j}\right) \left(2\,N +1\right)}}}{\textcolor{blue}{\sin \! \left(\frac{a}{2\,N +1}\right)}\textcolor{blue}{ \left(2\,N +1\right)}\,{\textcolor{red}{\sin \! \left(\frac{a}{\left(2\,N +1\right)^{2}}\right)\textcolor{red}{ \left(2\,N +1\right)}}}}=\frac{\sin \! \left(a \right)}{a} 
\label{G2xA}
\end{equation}
\fi

\begin{equation}
\textcolor{blue}{\left[\frac{{\sin \! \left(a \right)}\;}{\textcolor{blue}{\sin \! \left(\frac{a}{(2\,N +1)}\right)\left(2\,N +1\right)}}\right]}
\textcolor{red}{\left[\frac{\textcolor{red}{\sin \! \left(\frac{a}{2\,N\, +1}\right)}\;}{\textcolor{red}{\sin \! \left(\frac{a}{(2\,N +1)^2}\right)\left(2\,N +1\right)}}\right]}
{\underset{j =2}{\prod}}\textcolor{black}{ \frac{\sin \left(a/ \left(2\,N\, +1\right)^{j}\right)}{~\sin \left(a/ \left(2\,N +1\right)^{1+j}\right) \left(2\,N +1\right)}}=\frac{\sin \! \left(a \right)}{a} 
\label{G2xA}
\end{equation}

where the first two terms of the product have been extracted (and coloured) to demonstrate that the product partially telescopes (e.g.  the denominator of the j=0 factor (blue) cancels the corresponding factor in the numerator of the j=1 factor (red), each contributing a factor $2N+1$ in the denominator. The limiting term(s) of the product then approach unity because
\begin{equation}
\underset{j \rightarrow \infty}{\mathrm{lim}}\; \frac{1}{\sin \! \left(a/\left(2\,N +1\right)^{1+j}\right)}
 = 
\frac{ \left(2\,N +1\right)^{j+1}}{a}+O\! \left(\frac{1}{\left(2\,N +1\right)^{j}}\right)\,,
\label{LimXA}
\end{equation}
cancelling the accumulating multiplicative factor, thereby demonstrating the veracity of \eqref{G2x}, and similarly \eqref{G2A}.

\section{Some Perspectives} \label{sec:Perspective}

Infinite products have significance rooted in the early history of analysis exemplified by Euler's 1734 product formula \cite[Eq. 4.22.1]{NIST} for the sin function

\begin{equation}
\frac{\sin \! \left(a \right)}{a} = 
\overset{\infty}{\underset{j =1}{\prod}}\; \left(1-\frac{a^{2}}{\pi^{2}\,j^{2}}\right)\,.
\label{EuProd}
\end{equation}

In the realm of infinite products involving trigonometric functions, we have Euler's infinite product formula \eqref{GN1}. Setting $a=\pi/2$ in \eqref{GN1} yields
\begin{equation}
\frac{2}{\pi} = 
\frac{\sqrt{2}}{2}\times\frac{\sqrt{2+\sqrt{2}}}{2}\times\frac{\sqrt{2+\sqrt{2+\sqrt{2}}}}{2}\times\dots\,
\label{Viete}
\end{equation}
as a special case, originally obtained recursively by François Vi\`ete in 1593 and recently generalized by both Moreno and Garcia \cite{MorGar} and Levin \cite{Levin2005}, the latter of whose approach was based upon an an analysis of the functional equation

\begin{equation}
F \! \left(a\,z \right) = g \! \left(F \! \left(z \right)\right)
\label{Levin}
\end{equation}
with simple conditions on $F(z)~\mathrm{and}~g(z)$.\newline

Relevant to the identity reported here, in 1876, Dobinski \cite[Eq. (72)]{InfProd} proposed the identity
\begin{equation}
\overset{\infty}{\underset{j =1}{\prod}}\; \left(\tan\left(2^{j}\,a \right)\right)^{2^{-j}}
 = 4 \sin^{2}\left(a \right)\,.
\label{Dob1}
\end{equation}
Written otherwise, we have
\begin{equation}
\tan \! \left(a \right)\times \left({\tan}\left(2\,a \right)\right)^{\frac{1}{2}}\times \left(\tan\left(4\,a \right)\right)^{\frac{1}{4}}\times \left(\tan\left(8\,a \right)\right)^{\frac{1}{8}}\times\dots
 = 4 \sin^{2}\left(a \right)
\label{Dob1a}
\end{equation}
to be compared to \eqref{X1b} and \eqref{X1a}. In 1947, the identity \eqref{Dob1} was questioned by Agnew and Walker \cite{AgneWalk}. By reproducing Dobinski's derivation, Agnew and Walker claim that \eqref{Dob1} is true {\bf iff}

\begin{equation}
\underset{j \rightarrow \infty}{\mathrm{lim}}\! \left(\sin\left(2^{1+j}\,a \right)\right)^{2^{-j}}
 = 1\,;
\label{Agwa}
\end{equation} 
therefore, \eqref{Dob1} is true only for special values of $a\neq (2n+1)\pi/2^k$. To reinforce their scepticism, it is worth noting that for certain real values of $a$, some of the product factors are complex while the right-hand side is real. So, \eqref{Dob1} thus requires that, if for some value $j=j_{+}$ $\tan(2^{j_{+}}a)>0$, then $2^{j}a$ must lie in a region where $\tan(2^{j}a)>0$ for all $j \neq j_{+}$. As might be expected, this imposes a severe restriction on allowable values of $a$, and omits the possibility that the imaginary part of all the (complex) product factors corresponding to $\tan(2^{j_{+}}a)<0$ somehow vanishes. Unfortunately, Agnew and Walker do not identify a concrete example to demonstrate the limited veracity of \eqref{Dob1} and numerical experimentation by myself has uncovered no examples. 

\subsection{A digression}\label{sec:Digression}

The juxtaposition of the proof of \eqref{Dob1} in proximity with \eqref{BR114} leads the curious analyst to consider the question of the evaluation of \eqref{BR114} when its right-hand side vanishes, in the light of ``Weierstrass' Factor-Theorem", \cite[page 9, Theorem 1]{Knopp2} that states ``A convergent product has the value zero if and only if, one of its factors vanishes". Which of the factors of \eqref{BR114} vanishes in this eventuality? Consider the candidate case $a=2^n\pi$. After evaluating the appropriate limit on the right-hand side, we find

\begin{equation}
\overset{k -1}{\underset{j =0}{\prod}}\; \cos \! \left(2^{j+n}\,\pi \right)
 = 
1+\left(-\frac{2^{2\,k}}{6}+\frac{1}{6}\right) \left(a -2^{n}\,\pi \right)^{2}+O\! \left(\left(a -2^{n}\,\pi \right)^{4}\right)
\label{case1}
\end{equation}
and there is no issue, since each of the factors comprising the product are equal to unity. In the case $a=\pi/2^m$, that is

\begin{equation}
\overset{k -1}{\underset{j =0}{\prod}}\; \cos \! \left(\frac{2^{j}\,\pi}{2^{m}}\right)
 = 
\frac{\sin \! \left({2^{k-m}\,\pi}\right)}{2^{k}\,\sin \! \left({2^{-m}\pi}\right)}
\label{case2}
\end{equation}

only if $k\geq m$ does the right-hand side vanish, and accordingly, only the single left-hand side factor indexed by $j=m-1$ vanishes, and  both \eqref{BR114} and Weierstrass' Factor-Theorem are satisfied. Now consider \eqref{BR114} in the case that $k\rightarrow\infty$, that is

\begin{equation}
\underset{k \rightarrow \infty}{\mathrm{lim}}\; \overset{k -1}{\underset{j =0}{\prod}}\; \cos \! \left(2^{j}\,a \right)
 = 
\underset{k \rightarrow \infty}{\mathrm{lim}}\; \frac{\sin \! \left(2^{k}\,a \right)}{2^{k}\,\sin \! \left(a \right)}\rightarrow 0\,.
\label{case3}
\end{equation}
Unlike the previous two cases, \eqref{case3} is apparently analytically (and, for large values of $k$, numerically) valid for all values of $a$. Consequently, consistency of the left-hand side requires the existence of a solution to

\begin{equation}
2^{j-1}a/N=\pi
\label{caselInf}
\end{equation}
for general values of $a$, if $N$ is an odd integer. That \eqref{caselInf} is clearly impossible (e.g. $a$ is rational) suggests that \eqref{case3} provides a counterexample to Weierstrass' factor theorem -- an unlikely prospect. Is there a resolution?\newline 
\iffalse
First of all, we must ask if the infinite product \eqref{case3} converges. Denoting the general term of \eqref{case3} by $c_j(a)=\cos(2^j\,a)$, the test for convergence requires that there must exist some $j_{0}$ such that
\begin{equation}
|c_{j}(a).c_{j+1}(a)...c_{j+r}(a)-1|<\epsilon
\label{test}
\end{equation}
for all $j>j_{0}$, where $\epsilon>0$ and all $r\geq 1$. That this is impossible is easily seen by a counterexample. Choose $r=1$ and 
\fi

After due consideration, one must conclude that equating the right-hand side limit of \eqref{case3} to zero is unjustified (indicated by the use of the arrow), because the limit $k\rightarrow\infty$ yields an {\it accumulation}, rather than a {\it limit}, point, courtesy of the numerator factor that varies discontinuously for consecutive values of $k$ in a neighbourhood of zero acting as a limit. Similarly, the factors comprising the left-hand side of \eqref{case3} span an infinite number of possibilities, any one of which can approach arbitrarily close to, but not equal to, zero, thus emulating an accumulation point and the infinite product does not converge, although it comes arbitrarily close to zero. However, let us now rewrite \eqref{case3} as follows

\begin{equation}
\frac{\sin \! \left(a \right)\overset{k -1}{\underset{j =0}{\prod}}\; \cos \! \left(2^{j}\,a \right)}{\sin \! \left(2^{k}\,a \right)}
 = \frac{1}{2^{k}}\,
\label{BR114a}
\end{equation}
an identity that is also analytically and numerically valid. Again, considering the limit $k\rightarrow \infty$ gives
\begin{equation}
\underset{k \rightarrow \infty}{\mathrm{lim}}\;\frac{\sin \! \left(a \right)\overset{k -1}{\underset{j =0}{\prod}}\; \cos \! \left(2^{j}\,a \right)}{\sin \! \left(2^{k}\,a \right)}
 = \underset{k \rightarrow \infty}{\mathrm{lim}}\;\frac{1}{2^{k}}\,=0,
\label{CaseInfa}
\end{equation}
and this time there is no ambiguity about the limit of the right-hand side. However, as seen above, since the infinite product term on the left-hand side must be non-convergent, Weierstrass' factor theorem does not apply and we have an example of a simple equality whose left and right-hand side limits differ, since neither of the other two left-hand side factors will lead to a zero or negate the arbitrariness of the product term; if the product were to be convergent then we would have a counterexample to Weierstrass' factor theorem. Since $k$ is not a continuous variable, I am unaware that this example violates any theorems of analysis.\newline

In addition to the above, note that \eqref{Jo2} exhibits similar behaviour when $k\rightarrow \infty$, although in that case, the resolution is obvious since the right-hand side clearly alternates in value between even and odd values of $k$. However, \eqref{Jo2} does exhibit another curiosity: when $k$ is even, the left-hand side product possesses two vanishing factors indexed by $j=k/2$ and $j=3k/2$, whereas the right-hand side possess a simple zero. Again, since $k$ is not a continuous variable, the disagreement in the {\it order} of the zero (if one can talk about the order of a zero attached to a discrete variable) does not violate Weierstrass' theorem as written. For a collection of other interesting analytic anomalies, see \cite{CounterX}.

\section{Summary}

What began as a simple derivation of what appears to be a new, curious identity, led to a study of some related and well-known identities. This in turn led to a unification of heretofore disparate listed results for a family of infinite product identities characterized by arguments containing the index variable appearing as an exponent. Further investigation of related identities selected from the literature led to some interesting examples of unexpected analytic behaviour.

\begin{appendices}

\section{Appendix: Lemma} \label{sec:Appendix}

% From file ProdLimGen.mw and testFunctionalEq.mw
%{\bf Lemma:}
\begin{Lemma}
Consider the functional equation
\begin{equation}
f \! \left(a \right) = 
g \! \left(a \right)+x\,f \! \left(\frac{a}{p}\right)
\label{R0}
\end{equation}
Then, if $|x|< 1$, $f(0)$ is finite, $1<p\in\Re$ and if the series converges 
\begin{equation}
f \! \left(a \right) = 
\overset{\infty}{\underset{j =0}{\sum}}\; x^{j}\,g \! \left(\frac{a}{p^{j}}\right)
\label{R0000}
\end{equation}
{\bf Proof:}
Let $a:=a/p$ in \eqref{R0} and substitute into itself, to yield
\begin{equation}
f \! \left(a \right) = g \! \left(a \right)+g \! \left(\frac{a}{p}\right)\,x 
+ f \! \left(\frac{a}{\,p^2}\right)\,x^{2}\,.
\label{R00}
\end{equation}
Let $a:=a/p^2$ in \eqref{R0}, substitute into \eqref{R00} and repeat recursively altogether $N$ times to obtain
\begin{equation}
f \! \left(a \right) = 
 f \! \left(\frac{a}{p^{N}}\right) x^{N}+\overset{N-1}{\underset{j =0}{\sum}}\; x^{j}\,g \! \left(\frac{a}{p^{j}}\right)\,.
\label{R0N}
\end{equation}
According to the suppositions, if $N\rightarrow\infty$, the first term in \eqref{R0N} vanishes, leaving \eqref{R0000}. {\bf QED} 
\end{Lemma}
% from file TestFunctional Equation.mw

Corollary 1: If $x=1$ and $f(0)=0$, \eqref{R0000} also follows.

Corollary 2: If $x>1$, $\underset{x \rightarrow \infty}{\mathrm{lim}}\;f \left(x \right)$ is finite, and if the sum converges, by reversing \eqref{R0} to read
\begin{equation}
f \! \left(\frac{a}{p}\right) = 
-\frac{g \! \left(a \right)}{x}+\frac{f \! \left(a \right)}{x}
\label{L1}
\end{equation}

we have the solution
\begin{equation}
f \! \left(a \right) = 
-\overset{\infty}{\underset{j =1}{\sum}}\; \frac{g \! \left(p^{j}\,a \right)}{x^{j}}\,.
\label{LN1}
\end{equation}
\end {appendices}

\bibliographystyle{unsrt}

\bibliography{c://physics//biblio}

\end{flushleft}
\end{document}